\documentclass[10pt]{amsart}
\usepackage{latexsym, amsmath,amssymb}

\usepackage{epsfig}

\renewcommand{\epsilon}{\varepsilon}
\newcommand{\newsection}[1]
{\subsection{#1}\setcounter{theorem}{0} \setcounter{equation}{0}
\par\noindent}

\newtheorem{theorem}{Theorem}

\newtheorem{lemma}[theorem]{Lemma}
\newtheorem{corr}[theorem]{Corollary}

\newtheorem{proposition}[theorem]{Proposition}
\newtheorem{deff}[theorem]{Definition}

\newcommand{\bth}{\begin{theorem}}
\newcommand{\ble}{\begin{lemma}}
\newcommand{\bcor}{\begin{corr}}

\newcommand{\bdeff}{\begin{deff}}

\newcommand{\bprop}{\begin{proposition}}
\newcommand{\ele}{\end{lemma}}
\newcommand{\ecor}{\end{corr}}
\newcommand{\edeff}{\end{deff}}

\newcommand{\eprop}{\end{proposition}}

\renewcommand{\Pi}{\varPi}

\renewcommand{\epsilon}{\varepsilon}

\newcommand{\tidle}{\tilde}

\begin{document}

\title[Kakeya-Nikodym averages and  $L^p$-norms of eigenfunctions]
{Kakeya-Nikodym averages and \\ $L^p$-norms of eigenfunctions}
\thanks{The author was supported in part by NSF
Grant DMS-0555162.
}


\author{Christopher D. Sogge}
\address{Department of Mathematics,  Johns Hopkins University,
Baltimore, MD 21218}
\email{sogge@jhu.edu}

\maketitle

\begin{abstract}  We provide a necessary and sufficient condition that $L^p$-norms, $2<p<6$,
of eigenfunctions of the square root of minus the Laplacian on two-dimensional compact boundaryless Riemannian manifolds $M$  are small compared to a natural power of the eigenvalue $\lambda$.  The condition that ensures this is that their $L^2$-norms over $O(\lambda^{-1/2})$ neighborhoods of arbitrary unit geodesics are small when $\lambda$ is large (which is not the case for the highest weight spherical harmonics on $S^2$ for instance).  The proof exploits Gauss' lemma and the fact that the bilinear oscillatory integrals in 
H\"ormander's proof of the Carleson-Sj\"olin theorem become better and better behaved away from the diagonal.  Our results are related to a recent work of Bourgain who showed that $L^2$-averages over geodesics of eigenfunctions are small compared to a natural power of the eigenvalue $\lambda$ provided that the $L^4(M)$ norms are similarly small.  Our results imply that QUE cannot hold on a compact boundaryless Riemannian manifold $(M,g)$ of dimension two if $L^p$-norms are saturated for a given $2<p<6$.  We also show that eigenfunctions cannot have a maximal rate of $L^2$-mass concentrating along unit portions of geodesics that are not smoothly closed.
\end{abstract}

\newsection{Introduction}

The main purpose of this paper is to slightly sharpen a recent result of Bourgain~\cite{bourgainef}
concerning two-dimensional compact boundaryless Riemannian manifolds.  By doing so
we shall be able to provide a natural necessary and sufficient condition concerning the
growth rate of $L^p$-norms of eigenfunctions for $2<p<6$ and their $L^2$-concentration
about geodesics.

There are different ways of measuring the concentration of eigenfunctions.  One is by means
of the size of their $L^p$-norms for various values of $p>2$.  If $M$ is a compact
boundaryless manifold  with Riemannian metric $g=g_{jk}(x)$ and if
$\Delta_g$ is the associated Laplace-Beltrami operator, then the eigenfunctions
solve the equation $-\Delta_g e_{\lambda_j}(x)=\lambda^2_j e_{\lambda_j}(x)$ for
a sequence of eigenvalues $0=\lambda_0\le \lambda_1\le \lambda_2\dots$.  Thus, we are
normalizing things so that $\lambda_j$ are the eigenvalues of the first-order operator
$\sqrt{-\Delta_g}$.  We shall also usually assume that the $e_{\lambda_j}$ have $L^2$-norm
one, in which case $\{e_{\lambda_j}\}$ provides an orthonormal basis of $L^2(M,dx)$ where
$dx$ is the volume element coming from the metric.  Earlier, in the two-dimensional case,  we showed in \cite{soggeest} that if $M$ is fixed
then there is a uniform constant $C$ so that for $2\le p\le \infty$ and $j=1,2,3,\dots$ 
\begin{equation}\label{1}
\|e_{\lambda_j}\|_{L^p(M)}\le C\lambda_j^{\delta(p)}\|e_{\lambda_j}\|_{L^2(M)},
\end{equation}
with
$$\delta(p)= \begin{cases} \displaystyle
\frac12(\frac12-\frac1p), \quad 2\le p\le 6,
\\ \\
\displaystyle
\frac12 -\frac2p, \quad 6\le p\le \infty.
\end{cases}
$$

These estimates are sharp for the round sphere $S^2$, and in this case they detect
two types of concentration of eigenfunctions that occur there.  Recall that on $S^2$ with
the canonical metric the distinct eigenvalues are $\sqrt{k^2+k}$, $k=0,1,2,\dots$, which 
repeat with multiplicity $d_k=2k+1$.  If ${\mathcal H}_k$, the space of spherical
harmonics of degree $k$, is the space of all eigenfunctions with eigenvalue
$\sqrt{k^2+k}$, and if $H_k(x,y)$ is the kernel of the projection operator onto ${\mathcal H}_k$,
then the $k$-th zonal function at $x_0\in S^2$ is $Z_k(y)=(H_k(x_0,x_0))^{-1/2}H_k(x_0,y)$.
Its $L^2$-norm is one but its mass is highly concentrated at $\pm x_0$ where it takes
on the value $\sqrt{d_k/4\pi}$.  Explicit calculations show that $\|Z_k\|_{L^p(S^2)}\approx
k^{\delta(p)}$ for $p\ge 6$ (see e.g. \cite{sph}), which shows that in the case of
$M=S^2$ with the round metric \eqref{1} cannot be improved for this range of 
exponents.  Another extreme type of concentration is provided by the
highest weight spherical harmonics which have mass concentrated on the equators of
$S^2$, which are its geodesics.  The ones concentrated on the equator $\gamma_0=\{(x_1,x_2,0); 
\, x_1^2+x_2^2=1\}$ are the functions $Q_k$, which are the restrictions of the ${\mathbb R}^3$
harmonic polynomials $k^{1/4}(x_1+ix_2)^k$ to $S^2=\{x; |x|=1\}$.  One can check
that the $Q_k$ have $L^2$-norms comparable to one and $L^p$-norms
comparable to $k^{\frac12(\frac12-\frac1p)}$ when $2\le p\le 6$ (see e.g. \cite{sph}).
Notice also that the $Q_k$ have Gaussian type concentration about the equator
$\gamma_0$.  Specifically, if ${\mathcal T}_{k^{-1/2}}(\gamma_0)$ denotes all
points on $S^2$ of distance smaller than $k^{-1/2}$ from $\gamma_0$ then one can
check that
\begin{equation}\label{i.3}
\liminf_{k\to\infty}\int_{{\mathcal T}_{k^{-1/2}}(\gamma_0)} |Q_k(x)|^2 \, dx >0.
\end{equation}
For future reference, obviously the $Q_k$ also have the related property that
\begin{equation}\label{i.4}
\int_{\gamma_0}|Q_k|^2 \, ds \approx k^{1/2},
\end{equation}
 if $ds$ is the measure on $\gamma_0$ induced by the 
volume element.  

Thus, the sequence of highest weight spherical harmonics shows that the norms
in \eqref{1} (for $2<p<6$), \eqref{i.3} and \eqref{i.4} are related.  A goal of this
paper is to show that this is true for general two-dimensional compact manifolds
without boundary.

We remark that, although the estimates \eqref{1} are sharp for the round sphere,
one expects that  it should be the case that, for generic
manifolds, and $L^2$-normalized eigenfunctions one has
\begin{equation}\label{i.5}
\limsup_{j\to \infty} \lambda^{-\delta(p)}_j \|e_{\lambda_j}\|_{L^p(M)}=0
\end{equation}
for every $2<p\le \infty$.
  This was verified for exponents $p>6$ by Zelditch and the author in \cite{soggezelditch} by showing that if there are no points $x$  through which a positive
measure of geodesics starting at $x$ loop back through $x$ then $\|e_\lambda\|_\infty
=o(\lambda^{1/2})$.  By interpolating with the estimate \eqref{1} for $p=6$, this yields
\eqref{i.5} for all $p>6$.  Corresponding results were also obtained in \cite{soggezelditch}
for higher dimensions.  Recently, these results were strengthened by Toth, Zelditch and
the author \cite{stz} to allow similar results for quasimodes under the weaker condition that
at every point $x$ the set of recurrent directions for the first return map for geodesic flow has measure
zero in the cosphere bundle $S^*_xM$ over $x$.

Other than the partial results in Bourgain \cite{bourgainef}, there do not seem to be
any results addressing when \eqref{i.5} holds for a given $2<p<6$ (although Zygmund \cite{zygmund} showed that on the torus $L^2$-normalized eigenfunctions have
uniformly bounded
$L^4$-norms).  Furthermore, there do not seem to be results addressing the interesting
endpoint case of $p=6$, where one expects both types of concentration mentioned before to
be relevant.

Recently authors have studied the $L^2$ norms of eigenfunctions over unit-length geodesics.
Burq, G\'erard and Tzvetkov \cite{burq} showed that if $\Pi$ is the collection of all
unit
length geodesics then
\begin{equation}\label{i.50}\sup_{\gamma\in \Pi}\int_{\gamma}|e_{\lambda_j}|^2 \, ds
\lesssim \lambda^{1/2}_j\|e_{\lambda_j}\|^2_{L^2(M)}, \, \, j=1,2,3,\dots,
\end{equation}
which is sharp in view of \eqref{i.4}.  Related results for hyperbolic
surfaces were obtained earlier by Reznikov \cite{rez}, who opened up the present line of investigation.
The proof of \eqref{i.50} boils down to bounds for certain Fourier integral operators
with folding singularities (cf. Greenleaf and Seeger \cite{greenleafseeger}, Tataru \cite{tataru}).
In \S 3, we shall use ideas from \cite{greenleafseeger}, \cite{tataru}, and \cite{dg}, \cite{ivrii}, \cite{stz},
\cite{soggezelditch} to show that if $\gamma\in \Pi$ and
$$\limsup_{j\to\infty}\lambda_j^{-1/2}\int_\gamma |e_{\lambda_j}|^2 \, ds >0,$$
then the geodesic extension of $\gamma$ must be a smoothly closed geodesic.
Presumably it also has to be stable, but we cannot prove this.
Further recent work on $L^2$-concentration along curves can be found in Toth \cite{toth}.

In a recent paper \cite{bourgainef}, Bourgain proved an estimate that partially links the
norms in \eqref{1} and \eqref{i.50}, namely that for all $p\ge2$
\begin{equation}\label{bourgain}
\sup_{\gamma\in \Pi}\int_\gamma |e_{\lambda_j}|^2\, ds \lesssim \lambda^{1/p}_j
\|e_{\lambda_j}\|_{L^p(M)}^2.
\end{equation}
Of course for $p=2$, this is just \eqref{i.50}; however, an interesting feature of \eqref{bourgain}
is that the estimate for a given $2<p\le 6$ combined with \eqref{1} yields \eqref{i.50}.
Thus, if $e_{\lambda_{j_k}}$ is a sequence of eigenfunctions with (relatively) small $L^p(M)$ norms
for a given $2<p\le 6$, it follows that its $L^2$-norms over unit geodesics must also
be (relatively) small.  Bourgain \cite{bourgainef} also came close to establishing
the equivalence of these two things by showing that given $\varepsilon>0$ there
is a constant $C_\varepsilon$ so that for $j=1,2,\dots$
\begin{equation}\label{i.7}
\|e_{\lambda_j}\|_{L^4(M)}\le C_\varepsilon \Bigl(\, \lambda_j^{1/8+\varepsilon}
\|e_{\lambda_j}\|_{L^2(M)}\, \Bigr)^{3/4} \, \Bigl[ \, \lambda^{-1/2}_j
\sup_{\gamma\in \Pi} \int_\gamma |e_{\lambda_j}|^2\, ds \, 
\Bigr]^{1/8}.
\end{equation}
Since $\delta(4)=1/8$ in \eqref{1}, if the preceding inequality held for $\varepsilon=0$
one would obtain the linkage of the size of the norms in \eqref{i.50} for large energy
with the size of the $L^4(M)$ norms.  Our main estimate in Theorem~\ref{theorem1} 
is that a variant of \eqref{i.7} holds, which is strong enough to complete the linkage.

Bourgain's approach in proving \eqref{i.7} was to employ ideas going back to 
C\'ordoba \cite{cordoba} and Fefferman \cite{fefferman} that were used to give a proof
of the Carleson-Sj\"olin theorem \cite{carsj}.  The key object that arose in C\'ordoba's
work \cite{cordoba} was what he called the Kakeya maximal function in ${\mathbb R}^2$,
namely,
\begin{equation}\label{i.8}
{\mathcal M}f(x)=\sup_{x\in{\mathcal T}_{\lambda^{-1/2}}}|{\mathcal T}_{\lambda^{-1/2}}|^{-1}\int_{{\mathcal T}_{\lambda^{-1/2}}}
|f(y)|\, dy, \quad f\in L^2({\mathbb R}^2),
\end{equation}
with the supremum taken over all $\lambda^{-1/2}$-neighborhoods ${\mathcal T}_{\lambda^{-1/2}}$ of unit line segments containing $x$, and $|{\mathcal T}_{\lambda^{-1/2}}|
\approx \lambda^{-1/2}$ denoting its area.  The above maximal operator is now more 
commonly called the Nikodym maximal operator as this is the terminology in
Bourgain's important papers \cite{bourgain1}--\cite{bourgain3} which established  highly nontrivial progress towards establishing the higher dimensional version of 
the Carleson-Sj\"olin theorem for Euclidean spaces ${\mathbb R}^n$, $n\ge3$.

One could also consider variable coefficient versions of the maximal operators in
\eqref{i.8}.  In the present context if $\gamma\in \Pi$ is a unit geodesic, one could 
consider the $\lambda^{-1/2}$-tube about it given by
$${\mathcal T}_{\lambda^{-1/2}}(\gamma)=\{\, y\in M; \, 
\inf_{x\in \gamma} d_g(x,y)<\lambda^{-1/2}\},$$
with $d_g(x,y)$ being the geodesic distance between $x$ and $y$.  Then if
$\text{Vol}_g({\mathcal T}_{\lambda^{-1/2}}(\gamma))$ denotes the measure of this tube,
the analog of \eqref{i.8} would be
$${\mathcal M}f(x)=\sup_{x\in \gamma\in \Pi}
\frac1{\text{Vol}_g({\mathcal T}_{\lambda^{-1/2}}(\gamma))}\int_{{\mathcal T}_{\lambda^{-1/2}}}
|f(y)|\, dy.$$
These operators have been studied before because of their applications in harmonic
analysis on manifolds.  See e.g. \cite{mss}, \cite{sonick}.  As was shown in \cite{mins},
following the earlier paper \cite{bourgain3}, they are much better behaved in
2-dimensions compared to higher dimensions.

As \eqref{i.7} suggests, it is not the size of the $L^2$-norm of ${\mathcal M}f$ for
$f\in L^2(M)$ that is relevant for estimating $L^4(M)$-norms of eigenfunctions
but rather the sup-norm of this quantity with $f=|e_{\lambda_j}|^2$, which
up to the normalizing factor in front of the integral is the quanitity
$$\sup_{\gamma\in \Pi}\int_{{\mathcal T}_{\lambda^{-1/2}}(\gamma)}|e_{\lambda_j}(x)|^2\, dx.$$
If the $e_{\lambda_j}$ are $L^2$-normalized this is trivially bounded by one.  In rough terms
our results say that beating this trivial bound is equivalent to beating the bounds in
\eqref{1} for a given $2<p<6$.  

Let us now state our variant of \eqref{i.7}:

\begin{theorem}\label{theorem1}  Fix a two-dimensional compact boundaryless 
Riemannian manifold $(M,g)$.
Then given $\varepsilon>0$ there is a constant $C_\varepsilon$ so that for eigenfunctions $e_\lambda$
of $\sqrt{-\Delta_g}$ with eigenvalues $\lambda\ge1$ we have
\begin{multline}\label{t1}
\|e_\lambda\|_{L^4(M)}^4 \le \varepsilon\lambda^{1/2}\|e_\lambda\|_{L^2(M)}^4+
C_\varepsilon \lambda^{1/2}\|e_\lambda\|_{L^2(M)}^2\sup_{\gamma\in \Pi}
 \int_{{\mathcal T}_{\lambda^{-1/2}(\gamma)}}|e_\lambda(x)|^2\, dx 
\\
+C\|e_\lambda\|_{L^2(M)}^4,
\end{multline}
with $C$ being a fixed constant which is independent of $\lambda$ and $\varepsilon$.
\end{theorem}

We shall prove this not by adapting C\'ordoba's \cite{cordoba} proof of the Carleson-Sj\"olin theorem but rather that of H\"ormander \cite{hormander}.  He obtained sharp oscillatory
integral bounds in ${\mathbb R}^2$ that provided sharp B\"ochner-Riesz estimates for
$L^4({\mathbb R}^2)$ (i.e. the Carleson-Sj\"olin theorem), which turns out to be the endpoint
case for this problem in 2-dimensions.  H\"ormander's approach was to turn this $L^4$-problem
into an $L^2$-problem by squaring the oscillatory integrals and then estimating their $L^2$-norms.  As his proof shows, the resulting bilinear operators that arise are better and
better behaved away from the diagonal, and this fact is what allows us to take the constant
in front of the first term in the right side of \eqref{t1} to be arbitrarily small (at the expense
of the 2nd term).

Stein~\cite{stein} provided a generalization of H\"ormander's oscillatory integral theorem
to higher dimensions in a way that proved to be sharp because of a later construction
of Bourgain~\cite{bourgain3}.  Bourgain's example and related ones in \cite{mins} suggest
that extending the results of this paper to higher dimensions (where the range of exponents
would be $2<p<2(n+1)/(n-1)$) could be subtle.  On the other hand, since the constructions
tend to involve concentration about hypersurfaces as opposed to geodesics, their relevance
is not plain.

We shall prove Theorem~\ref{theorem1} by estimating an oscillatory integral operator, which up to a remainder
term, reproduces eigenfunctions.  The remainder term  in this reproducing formula accounts for the last term in \eqref{t1},
which we could actually take to be $\le C_N\lambda^{-N}\|e_\lambda\|_2^4$ for any $N$,
but this is not important for our applications.  Also, we remark that the proof of the Theorem
will show that the constant $C_\varepsilon$ in \eqref{t1} can be taken to be $O(\varepsilon^{-2})$
as $\varepsilon\to 0$.

Let us now state an immediate consequence of Theorem \ref{theorem1} which states
that the size of $L^4$-norms of eigenfunctions is equivalent to size of $L^2$-mass near
geodesics.

\begin{corr}\label{corollary1}  Let $e_{\lambda_{j_k}}$ be a sequence of eigenfunctions
with eigenvalues $\lambda_{j_1}\le \lambda_{j_2}\le \dots$ and unit $L^2(M)$-norms.  Then
\begin{equation}\label{i}
\limsup_{k\to\infty} \sup_{\gamma\in \Pi}\int_{{\mathcal T}_{\lambda^{-1/2}_{j_k}}(\gamma)}
|e_{\lambda_{j_k}}(x)|^2\, dx =0
\end{equation}
if and only if
\begin{equation}\label{ii}
\limsup_{k\to \infty}\lambda_{j_k}^{-1/8}\|e_{\lambda_{j_k}}\|_{L^4(M)}=0. \end{equation}
\end{corr}

To prove this, we first notice that if we assume \eqref{i}, then \eqref{ii} must hold because
of \eqref{t1}.  Also, by H\"older's inequality
$$\Bigl(\, \int_{{\mathcal T}_{\lambda^{-1/2}}(\gamma)}
|e_{\lambda}(x)|^2\, dx\, \Bigr)^{1/2}\le \bigl(\, \text{Vol}_g({\mathcal T}_{\lambda^{-1/2}}(\gamma))\,
\bigr)^{1/4}\|e_\lambda\|_{L^4(M)}\lesssim \lambda^{-1/8}\|e_\lambda\|_{L^4(M)},$$
and so \eqref{ii} trivially implies \eqref{i}.

If we use Bourgain's estimate \eqref{bourgain} and \eqref{1} we can say a bit more.

\begin{corr}\label{corollary2}  Let $\{e_{\lambda_{j_k}}\}_{k=1}^\infty$ be as above and suppose
that $2<p<6$.  Then the following are equivalent
\begin{align}
\label{bi}
\limsup_{k\to\infty}\lambda^{-1/2}_{j_k}\sup_{\gamma\in\Pi}\int_{\gamma}
|e_{\lambda_{j_k}}(s)|^2\, ds&=0
\\
\label{bii}
\limsup_{k\to\infty}\sup_{\gamma\in\Pi}\int_{{\mathcal T}_{\lambda^{-1/2}_{j_k}}(\gamma)}
|e_{\lambda_{j_k}}(x)|^2 \, dx &=0
\\
\label{biii}
\limsup_{k \to \infty}\lambda_{j_k}^{-\delta(p)}\|e_{\lambda_{j_k}}\|_{L^p(M)}&=0.
\end{align}
\end{corr}

To prove this result, we first note that, by the M. Riesz interpolation theorem and \eqref{1} for $p=2$ and $p=6$, \eqref{biii} holds for a given $2<p<6$ if and only if it holds for $p=4$, which
we just showed is equivalent to \eqref{bii}.  Clearly \eqref{bi} implies \eqref{bii}.  Finally,
since Bourgain's estimate \eqref{bourgain} shows that \eqref{biii} implies \eqref{bi},
the proof of Corollary~\ref{corollary2} is complete.

Let us conclude this section by describing one more application.  Recall that a sequence of $L^2$-normalized eigenfunctions $\{e_{\lambda_{j_k}}\}_{k=1}^\infty$ satisfies the quantum unique ergodicity property (QUE) if the associated Wigner measures $|e_{\lambda_{j_k}}|^2dx$ tend to the Liouville measure on $S^*M$.  If this is the case, then one certainly cannot have
$$\limsup_{k\to \infty} \sup_{\gamma\in\Pi}\int_{{\mathcal T}_{\lambda^{-1/2}_{j_k}}(\gamma)}
|e_{\lambda_{j_k}}(x)|^2 \, dx>0,$$
since the tubes are shrinking.  

In the case where $M$ has negative sectional curvature Schnirelman's \cite{schnirelman} theorem, proved by Zelditch \cite{zelditch}, says there is a density one subsequence $\{e_{\lambda_{j_k}}\}_{k=1}^\infty$
of all the $\{e_{\lambda_j}\}$ satisfying QUE.  Rudnick and Sarnak \cite{sarnak} conjectured that in the negatively curved
case there should be no exceptional subsequences violating QUE, i.e., in this case QUE should hold for the full sequence $\{e_{\lambda_j}\}$ of $L^2$-normalized eigenfunctions.  On
the other hand, by Corollary \ref{corollary2}, we have the following.

\begin{corr}\label{corollary4}  Let $M$ be a two-dimensional compact boundaryless
Riemannian manifold.  Then QUE cannot hold for $M$ if for a given $2<p<6$ there is saturation of $L^p$ norms, i.e.
$$\limsup_{j\to\infty}\lambda_j^{-\delta(p)}\|e_{\lambda_j}\|_{L^p(M)}>0,$$
with $e_{\lambda_j}$ being the $L^2$-normalized eigenfunctions.
\end{corr}

See e.g. \cite{zelind} for connections between QUE and the Lindel\"of hypothesis, and see \cite{cdv} for recent developments regarding the QUE conjecture.

\newsection{Proof of Theorem 1: Gauss' lemma and the Carleson-Sj\"olin condition}

As in \cite{bourgainef}  and  \cite{burq} we shall prove our estimate by using
certain convenient operators that reproduce eigenfunctions.  Specifically,
we shall use a slight variant of a result from \cite{soggebook}, Chapter 5 that
was presented in \cite{burq}.

\begin{lemma}\label{reproduce}  Let $\delta>0$ be smaller than half of the injectivity radius
of $(M,g)$.  Then there is a function $\chi\in{\mathcal S}({\mathbb R})$ with $\chi(0)=1$
so that if $d_g(x,y)$ is the geodesic distance between $x,y\in M$
\begin{equation}\label{r1}
\chi_\lambda f(x)=\chi\bigl(\sqrt{-\Delta_g}-\lambda\bigr)f(x)
=\lambda^{1/2}\int_M e^{i\lambda d_g(x,y)}\alpha(x,y,\lambda) f(y)\, dy+R_\lambda f(x),
\end{equation}
where
$$
\|R_\lambda f\|_{L^\infty(M)}\le C_N \lambda^{-N}\|f\|_{L^1(M)}, \, \, \text{for all} \, \,  N=1,2,\dots,
$$
and $\alpha \in C^\infty$ has the property that
$$|\partial^\alpha_{x,y}\alpha(x,y,\lambda)|\le C_\alpha , \, \, \text{for all} \, \,  \alpha,$$
and, moreover,
\begin{equation}\label{r3}
\alpha(x,y,\lambda)=0 \, \, \text{if } \, d_g(x,y)\notin (\delta/2,\delta).
\end{equation}
\end{lemma}

Since $\chi_\lambda e_\lambda=e_\lambda$ and since the 4th power of the $L^4$-norm of
$R_\lambda e_\lambda$ is dominated by the last term in \eqref{t1}, we conclude that
in order to prove Theorem \ref{theorem1} it is enough to show that, given $\varepsilon>0$
there is a constant $C_\varepsilon$ so that when $\lambda\ge1$
\begin{multline}\label{t1'}
\int_M \left| \,\lambda^{1/2} \int_M e^{i\lambda d_g(x,y)}\alpha(x,y,\lambda) f(y)\, dy\, \right|^2 \,
|f(x)|^2 \, dx
\le \varepsilon \lambda^{1/4}\|f\|_{L^2(M)}^2\|f\|_{L^4(M)}^2
\\
+C_\varepsilon \lambda^{1/2}\|f\|_{L^2(M)}^2
\sup_{\gamma\in \Pi}\int_{{\mathcal T}_{\lambda^{-1/2}}(\gamma)}|f(x)|^2 \, dx,
\end{multline}
for, if $f=e_\lambda$, the first term in the right is bounded by a fixed constant times
$\varepsilon \lambda^{1/2}\|e_\lambda\|_{L^2(M)}^4$, because of \eqref{1}.

After applying a partition of unity (and abusing notation a bit), we may assume that
in addition to \eqref{r3}, $\alpha(x,y,\lambda)$ vanishes unless $x$ is in a small neighborhood of
some $x_0\in M$ and $y$ is in a small neighborhood of some $y_0\in M$ with $\delta/2<d_g(x_0,y_0)< \delta$.  We may assume both of these neighborhoods are contained in
the geodesic ball $B(x_0,10\delta)=\{y\in M; \, d_g(x_0,y)<10\delta\}$.  As mentioned before,
we are also at liberty to take $\delta>0$ to be small.

To simplify the calculations to follow, it is convenient to choose a natural coordinate system.  Specifically, we shall choose Fermi normal coordinates about the geodesic $\gamma_0$ which passes through $x_0$ and is perpendicular to the geodesic connecting $x_0$ and $y_0$.  These coordinates will be well defined on $B(x_0,10\delta)$ if $\delta$ is small.  Furthermore, we may assume that the image of $\gamma_0\cap B(x_0,10\delta)$ in the resulting coordinates is a line
segment which is parallel to the 2nd coordinate axis and that all horizontal line segments $s\to \{(s,t_0)\}$ are geodesic with the property that $d_g((s_1,t_0),(s_2,t_0))=|s_1-s_2|$.  See Figure 1 below.

\begin{figure*}[h]
\centering
\includegraphics[width=0.3\textwidth]{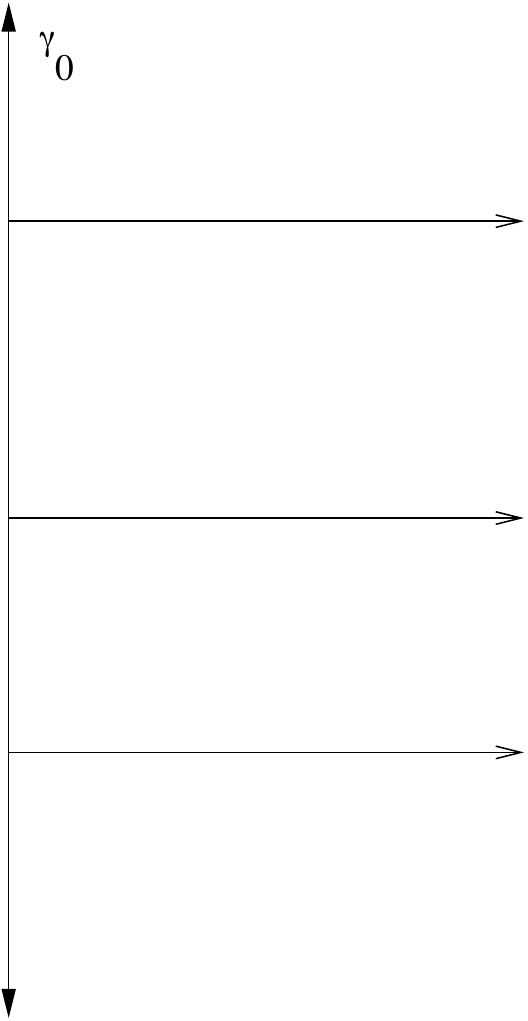}
\caption{Fermi normal coordinates about $\gamma_0$}
\end{figure*}

If we use these coordinates and apply Schwarz's inequality, we conclude that, in order to prove
\eqref{t1'}, it suffices to show that given $\varepsilon>0$ we can find $C_\varepsilon<\infty$ so that when $\lambda\ge1$
\begin{multline*}
\int\left(\, \int\Bigl|\, \lambda^{1/2}\int e^{i\lambda d_g(x,(s,t))}\alpha(x,(s,t),\lambda) f(s,t)\, dt\, \Bigr|^2 \, 
|f(x)|^2 \, dx \, \right) \, ds
\\
\le 
 \varepsilon \lambda^{1/4}\|f\|_{L^2(M)}^2\|f\|_{L^4(M)}^2
+C_\varepsilon \lambda^{1/2}\|f\|_{L^2(M)}^2
\sup_{\gamma\in \Pi}\int_{{\mathcal T}_{\lambda^{-1/2}}(\gamma)}|f(x)|^2 \, dx.
\end{multline*}
This, in turn would follow if we could show that given $\varepsilon>0$
\begin{multline}\label{2.5}
\int\, \Bigl| \lambda^{1/2}\int e^{i\lambda d_g(x,(s,t))}\alpha(x,(s,t),\lambda) \, h(t)\, dt \, \Bigr|^2 \, |f(x)|^2 \, dx
\\
\le \varepsilon \lambda^{1/4}\|h\|_{L^2(dt)}^2\|f\|_{L^4(M)}^2
+C_\varepsilon \lambda^{1/2}\|h\|^2_{L^2(dt)}\sup_{\gamma\in \Pi}\int_{{\mathcal T}_{\lambda^{-1/2}}(\gamma)} |f(x)|^2\, dx,
\end{multline}
with $C_\varepsilon$ depending on $\varepsilon>0$ but not on $s$ or on $\lambda\ge1$.

To simplify the notation, we shall establish this estimate for a particular value of $s$, 
which, after relabeling, we may assume to be $s=0$.  Since the proof of \eqref{2.5} for this case 
relies only on Gauss' lemma and the related Carleson-Sj\"olin condition, it also yields the
uniformity in $s$, assuming, as we may, that $\alpha$ has small support.

To prove this inequality, let us choose a function $\eta\in C^\infty_0({\mathbb R})$ satisfying $\eta(t)=0$, $|t|>1$, and $\sum_{j=-\infty}^\infty \eta(t-j)\equiv 1$.  
Given $\lambda\ge1$ fixed, we shall then set
$$\eta_j(t)=\eta_{\lambda,j}(t)=\eta(\lambda^{1/2}t-j).$$
Then, given $N=1,2,\dots$, we have that
\begin{multline}\label{N}
\left| \, \lambda^{1/2}\int e^{i\lambda d_g(x,(0,t))}\alpha(x,(0,t),\lambda)h(t) dt\, \right|^2
\\
\le N\sum_j\left|\lambda^{1/2}\int e^{i\lambda d_g(x,(0,t))}\eta_j(t)\alpha(x,(0,t),\lambda)h(t)dt
\right|^2
\\
+\left| \, \lambda \iint e^{i\lambda(d_g(x,(0,t))+d_g(x,(0,t'))}a_N(x,t,t')h(t)h(t')\, dt dt'\, \right|,
\end{multline}
where 
$$a_N(x,t,t')=\sum_{|j-k|>N}\eta_j(t)\alpha(x,(0,t),\lambda)\eta_k(t')\alpha(x,(0,t'),\lambda)$$
vanishes when $|t-t'|\le (N-1)\lambda^{-1/2}$.
The first term in the right side of the preceding inequality comes from applying Young's
inequality to handle the double-sum over indices with $|j-k|\le N$.  Because of \eqref{N}, we conclude
that \eqref{2.5} would follow if we could show that there is a constant
independent of $\lambda\ge 1$ and $N=2,3,4\dots$ so that
\begin{multline}\label{2.6}
\left\|\, \lambda \iint e^{i\lambda [d_g(x,(0,t))-d_g(x,(0,t')]}a_N(x,t,t')h(t)h(t') dt dt'\, \right\|_{L^2(dx)}
\\
\le C\lambda^{1/4}N^{-1/2}\|h\|_{L^2(dt)}^2,
\end{multline}
and also that there is a constant $C$ independent of $j\in {\mathbb Z}$ and $\lambda\ge1$ so that
\begin{multline}\label{2.7}
\int \Bigl| \, \lambda^{1/2}\int e^{i\lambda d_g(x,(0,t))}\eta_j(t) \, \alpha(x,(0,t),\lambda)h(t)\, dt \Bigr|^2
\, |f(x)|^2 \, dx
\\
\le C\lambda^{1/2}\|h\|^2_{L^2(dt)}\sup_{\gamma \in \Pi}\int_{{\mathcal T}_{\lambda^{-1/2}}(\gamma)}
|f(x)|^2 \, dx.
\end{multline}
Indeed, by using the finite overlapping of the supports of the $\eta_j$, if we set $\varepsilon=CN^{-1/2}$, then we see that these two inequalities and \eqref{N} imply  \eqref{2.5} with $C_\varepsilon \approx \varepsilon^{-2}$.  Since the proof of \eqref{2.7} only uses Gauss' lemma
and the fact that coordinates have been chosen so that $s\to(s,t_0)$ are unit speed geodesics for
fixed $t_0$, we shall just verify \eqref{2.7} for $j=0$, as the argument for this case will yield
the other cases as well.

The next step is to see that these two inequalities are consequences of the following two
propositions.

\begin{proposition}\label{csprop}  Let $a(x,t,t')$, $x\in {\mathbb R}^2$, $t,t'\in {\mathbb R}$
satisfy $|\partial^\alpha_x a|\le C_\alpha$ for all multi-indices $\alpha$ and
$a(x,t,t')=0$ if $|x|>\delta$ or $|t-t'|>\delta$ where $\delta>0$ is small.  Suppose also
that $\phi\in C^\infty({\mathbb R}^2\times{\mathbb R})$ is real and satisfies the
Carleson-Sj\"olin condition on the support of $a$, i.e.,
\begin{equation}\label{cs1}
\det \left( \begin{array}{cc}
\phi''_{x_1t} & \phi''_{x_2t} \\
\phi'''_{x_1tt} & \phi'''_{x_2tt} \\
 \end{array} \right) \ne 0.
 \end{equation}
 Then if the $\delta>0$ above is sufficiently small, there is a uniform constant $C$ so that
 when $\lambda, N\ge1$
 \begin{multline}\label{cs.2}
 \left\| \, \iint_{|t-t'|\ge N\lambda^{-1/2}}
 e^{i\lambda[\phi(x,t)+\phi(x,t')]}a(x,t,t')F(t,t') \, dt dt'\, 
 \right\|^2_{L^2({\mathbb R}^2)}
 \\
 \le C\lambda^{-3/2}N^{-1}\|F\|_{L^2({\mathbb R}^2)}^2.
 \end{multline}
 \end{proposition}
 
 To state the next Proposition, we need to introduce one more coordinate system, which 
 finally explains
where the $L^2$ norms over small tubular neighborhoods of geodesics comes into
play.  Since we are proving \eqref{2.7} with $j=0$ and since $\eta_0$ is supported 
in the small interval $[-\lambda^{-1/2},\lambda^{-1/2}]$, it is natural to take geodesic normal coordinates
about $(0,0)$.  If we recall that the 1st coordinate axis is a unit-speed geodesic in our original
Fermi normal coordinates, we shall naturally choose the geodesic normal coordinates $x\to
\kappa(x)$ that preserve this axis (and its orientation).  Such a system is unique up to reflection about this
axis, and we shall just fix one of these two choices.
 
 \begin{proposition}\label{gprop}  Let $\psi(x,t)=d_g\bigr(x, (0,t)\bigl)$, and suppose that $ \rho\in
C^\infty_0({\mathbb R}\times {\mathbb R}^2)$ satisfies
\begin{equation}\label{l.1}
|\partial^m_t\rho (t;x)|\le C_m( \lambda^{1/2})^m\, , \quad
\text{and}, \, \rho(t;x)=0, \, \, |t|\ge\lambda^{-1/2}.
\end{equation}
Suppose also that $\rho$ vanishes when $x$ is outside of a small neighborhood
${\mathcal N}$ of a fixed point $(-s_0,0)$ (in the Fermi normal coordinates) with $s_0>0$.
If $x\to \kappa(x)=(\kappa_1(x),\kappa_2(x))$ are the coordinates described above, assume
that points $x_j\in {\mathcal N}$ are chosen so 
that 
\begin{equation}\label{l.2}
\Bigl| \, \frac{\kappa_2(x_j)}{|\kappa(x_j)|} \, - \, \frac{\kappa_2(x_k)}{|\kappa(x_k)|} \, \Bigr|
\ge c \lambda^{-1/2}|j-k|, \quad \text{if } \, |j-k|\ge 10,
\end{equation}
with $c>0$ fixed.  It then follows that, if ${\mathcal N}$ is sufficiently small, then there is a uniform constant $C$, which is
independent of the $\{x_j\}$ chosen as above, so that
\begin{equation}\label{l.3}
\lambda^{1/2} \int \Bigl| \, \sum_j e^{i\lambda\psi(x_j,t)} \rho(t;x_j)\, a_j\, \Bigr|^2 \,  dt
\le C\sum |a_j|^2.
\end{equation}
\end{proposition}

Proposition \ref{csprop} would imply \eqref{2.6} if $\phi(x,t)=d_g(x,(0,t))$ satisfies the
Carleson-Sj\"olin condition.  The fact that this is the case is well known.  See e.g., Section 5.1 in \cite{soggebook}.  It follows from our choice of coordinates and the fact that if $x_0\in M$ is fixed then the set of points $\{\nabla_x d_g(x,y); \, x=x_0, \, d_g(x_0,y)\in (\delta/2,\delta)\}$ is the 
cosphere at $x_0$, $S^*_{x_0}M=\{\xi; \, \sum g^{jk}(x_0)\xi_j\xi_k=1\}$, where $g^{jk}(x)$
is the cometric (inverse to $g_{jk}(x)$).  If we choose geodesic normal coordinates $\kappa(y)$ vanishing at $x_0$ then the gradient becomes $\kappa(y)$.
This turns out to be equivalent to the usual formulation of Gauss' lemma, saying that this exponential map $y\to \kappa(y)$ is a local radial isometry.  More specifically, it says that small geodesic spheres centered at $x_0$ get
sent to spheres centered at the origin and small geodesic rays through $x_0$ intersect these geodesic spheres orthogonally and get sent to rays through the origin, which is what allows 
Proposition~\ref{gprop} to be true. 

Let us next see that Proposition \ref{gprop} implies \eqref{2.7} for $j=0$.  If we take
$\rho(t;x)=\eta_0(t)\alpha(x,(0,t),\lambda)$, then $\rho$ satisfies \eqref{l.1}.  Also, if we
let 
$$S_j=\{y; \, \theta(y)\in (\lambda^{-1/2}j,\lambda^{-1/2}(j+1)]\},$$
 where $\theta(y)\in [0,2\pi)$ is defined so that $y=|y|(\cos\theta(y), \sin\theta(y))$, then, if
 $y=\kappa(x)$ are the geodesic normal coordinates about $(0,0)$ in the Proposition~\ref{gprop}, then  the left side of \eqref{2.7} is dominated by
\begin{multline*}
\sum_j \Bigl\| \, \lambda^{1/2} \int e^{i\lambda \psi(x,t)}\rho(t;x)h(t) dt\Bigr\|_{L^\infty(\kappa^{-1}(S_j))}^2
\|f\|^2_{L^2(\kappa^{-1}(S_j)\cap K)}
\\
\le \sup_k \|f\|^2_{L^2(\kappa^{-1}(S_k)\cap K)} \sum_j \Bigl\|\lambda^{1/2}\int e^{i\lambda \psi(x,t)}\rho(t;x) h(t)dt\Bigr\|^2_{L^\infty(\kappa^{-1}(S_j))},
\end{multline*}
where $K$ is the $x$-support of $\rho$.  Since the first factor on the right is dominated by the
last factor in the right hand side of \eqref{2.7} (the sup can just be taken over $(0,0)\in \gamma\in \Pi$ here), we conclude that we would obtain this
inequality if we could show that there is a uniform constant so that for all choices of $x_j\in \kappa^{-1}(S_j)$
\begin{equation}\label{2.8}
\lambda^{1/2}\sum_j \left|\, \int e^{i\lambda\psi(x_j,t)}\rho(t;x_j)h(t) \, dt\, \right|^2 \le C\|h\|^2_{L^2(dt)}.
\end{equation}
This inequality is an estimate for an operator from $L^2(dt)\to \ell^2$.  The dual operator
is the one in Proposition \ref{gprop}.
Therefore since, by duality, \eqref{2.8} follows from \eqref{l.3} we get \eqref{2.7}.  To verify this assertion, we use
the fact that if $\rho$ has small support then the terms in \eqref{2.8} with $\rho(t;x_j)\neq 0$
will fulfill the hypotheses in Proposition \ref{gprop}.

To finish the proof of Theorem \ref{theorem1}  we must prove the two propositions.  Let
us start with the first one since it is pretty standard.  It is based on the well known fact that the 
bilinear oscillatory integrals arising in H\"ormander's \cite{hormander} proof of the Carleson-Sj\"olin \cite{carsj} theorem become better and better behaved away from the diagonal.

\medskip

\noindent{\bf Proof of Proposition \ref{csprop}:}  Let $\Phi(x;t,t')=\phi(x,t)+\phi(x,t')$ be the phase function in \eqref{cs.2}.  Then 
$\Phi$ is a symmetric function in the $(t,t')$ variables.  So if we make the change of variables
$$u=(t-t', t+t'),$$
then since $|du/d(t,t')|=2$, we see that \eqref{cs1} implies that the Hessian determinant of $\Phi$ satisfies
$$\left| \, \det \left( \, \frac{\partial^2\Phi}{\partial x\partial u}\right) \, \right| \ge c|u_1|,$$
for some $c>0$ on the support of $a$, if the latter is small.  Since
$\Phi(x;u)$ is an even function of the diagonal variable $u_1$, it must
be a $C^\infty$ function of $u_1^2$.  So if we make the final change of variables
$$v=\Bigl(\frac12 u^2_1,u_2\Bigr),$$
then since $|dv/du|=|u_1|$, it follows that $$\left| \, \det \left( \, \frac{\partial^2\Phi}{\partial x\partial v}\, \right)\right|\ge c,$$
for some $c>0$.
This in turn implies that if $v$ and $\tilde v$ are close then
$$\left| \, \nabla_x\bigl[ \, \Phi(x,v)-\Phi(x,\tilde v)\, \bigr]\right| \ge c'|v-\tilde v|,$$
for some $c'>0$,
and since $x,v\to \Phi$ is smooth, we also have that
$$\left| \, \partial_x^\alpha \bigl[ \, \Phi(x,v)-\Phi(x,\tilde v)\, \bigr]\right| \le C_\alpha
|v-\tilde v|,$$
for all multi-indices $\alpha$.  Therefore, if we let
$$K_\lambda(v,\tilde v)=\int_{{\mathbb R}^2}a(x,t,t')\overline{a(x,\tilde t,\tilde t')}
e^{i\lambda[\Phi(x,v)-\Phi(x,\tidle v)]} \, dx,$$
then by integrating by parts, we find that if the number $\delta>0$ in the statement
of the Proposition is small then for
$j=1,2,3,\dots$
\begin{align}\label{cs.3}
|K_\lambda(v,\tilde v)|&\le C_j(1+\lambda|v-\tilde v|)^{-2j}
\\
&\le C_j(1+\lambda |(t+t')-(\tilde t+\tilde t')|)^{-j}
(1+\lambda|(t-t')^2-(\tilde t-\tilde t')^2|)^{-j}.
\notag
\end{align}
Note that the left side of \eqref{cs.2} equals
$$\int\cdots\int_{|t-t'|, |\tilde t-\tilde t'|\ge N\lambda^{-1/2}}
K_\lambda(t,t';\tilde t,\tilde t')F(t,t')\overline{F(\tilde t,\tilde t')}dt dt' d\tilde t d\tilde t'.$$
We next claim that there is a uniform constant $C$ so that for $\lambda,N\ge1$
\begin{equation}\label{cs.4}
\sup_{\tilde t,\tilde t'}\int_{|t-t'|\ge N\lambda^{-1/2}} |K_\lambda| \,dtdt' \, , \, \, 
\sup_{t,t'}\int_{|\tilde t-\tilde t'|\ge N\lambda^{-1/2}} |K_\lambda| \, d\tilde t d\tilde t' \le C
\lambda^{-2}(\lambda^{1/2}/N).
\end{equation}
This follows from \eqref{cs.3} and the fact that if $\tau = s^2$ then $2sds=d\tau$ and so, given
$\tau_0\in {\mathbb R}$, we have
\begin{multline*}
\int_{s\ge N\lambda^{-1/2}}(1+\lambda |s^2-\tau_0|)^{-2}\, ds
=\frac12 \int_{\sqrt{\tau} \ge N\lambda^{-1/2}}(1+\lambda|\tau-\tau_0|)^{-2}\, \frac{d\tau}{\sqrt \tau}
\\
\le (\lambda^{1/2}/N)\int_{-\infty}^{+\infty}(1+\lambda|\tau|)^{-2}\, d\tau \le
C\lambda^{-1}(\lambda^{1/2}/N).
\end{multline*}
Since \eqref{cs.4} and Young's inequality yield \eqref{cs.2}, the proof is complete.  \qed

\medskip

To finish our task we need to prove the other Proposition, which is a straightforward application
of Gauss' lemma.

\medskip

\noindent {\bf Proof of Proposition \ref{gprop}:}
The support assumptions on the amplitude will allow us to linearize the function
$t\to \psi$ in the proof, which is a tremendous help.  Specifically,
$$\psi(x,t)=\psi(x,0)+ t(\partial_t\psi(x,0))+r(x,t),$$
where
\begin{equation}\label{l.4}
|\partial^m_t r(x,t)|\le C_m |t|^{2-m}, \, 0\le m\le 2,  \quad
\text{and } \, |\partial^m_tr|\le C_m, \quad m\ge2.
\end{equation}
Our choice of coordinates implies that 
$$\partial_t\psi(x,0)= \langle \nu, \kappa(x)/|\kappa(x)|\, \rangle,$$
where the inner-product is the euclidean one and $\nu \in {\mathbb R}^2$
is chosen so that
 $\langle \nu,\nabla\rangle$ is the pushforward of $\partial/\partial x_2$ at $(0,0)$ 
 under the map $x\to \kappa(x)$---i.e., tangent vector to the curve $t\to \kappa((0,t))$.  Since the pushforward of 
 $\partial/\partial x_1$ is itself under this map, it follows that the second coordinate of $\nu$ is 
nonzero.  (See Figure 2 below.)
Therefore, if ${\mathcal N}\ni (s_0,0)$
is small enough, then our assumption \eqref{l.2} implies that
\begin{equation}
|\partial_t\psi(x_j,0)-\partial_t\psi(x_k,0)|\ge c'\lambda^{-1/2}|j-k|, \quad
\text{if } |j-k|\ge10, \, \text{and } \, x_j,x_k\in {\mathcal N},
\end{equation}
for some constant $c'>0$.

It is easy now to finish the proof of \eqref{l.3}.  If we let
$$\rho(x_j,x_k;t)=\rho(t;x_j)\overline{\rho(t;x_k)}e^{i\lambda(\psi(x_j,0)+r(x_j,t))}
e^{-i\lambda(\psi(x_k,0)+r(x_k,t))},$$
it follows from \eqref{l.1} and \eqref{l.4} that
$$
|\partial^m_t\rho(x_j,x_k;t)|\le C_m\lambda^{m/2}, 
$$
and
$$ 
\rho(x_j,x_k;t)=0, 
\\ \text{if } |t|\ge \lambda^{-1/2}, \, x_j\notin {\mathcal N}, \, \text{or } \, 
x_k\notin {\mathcal N}.
$$
We can use this since the left side of \eqref{l.3} equals
$$\lambda^{1/2}\sum_{j,k}a_j\overline{a_k} 
\, \Bigl(\,  \int e^{i t \lambda(\partial_t\psi(x_j,0)-\partial_t\psi(x_k,0))}\rho(x_j,x_k;t)\, dt\, \Bigr),
$$
which, after integrating by parts $N=1,2,3\dots$ times, we conclude is dominated by
a fixed constant $C_N$ times
$$\sum_{j,k}|a_ja_k|\, \bigl(\, 1+|j-k|\, \bigr)^{-N}.$$
Since, by Young's inequality, this is dominated by the right side of \eqref{l.3} when $N=2$, the proof is complete. \qed

\begin{figure*}[h]
\centering
\includegraphics[width=0.75\textwidth]{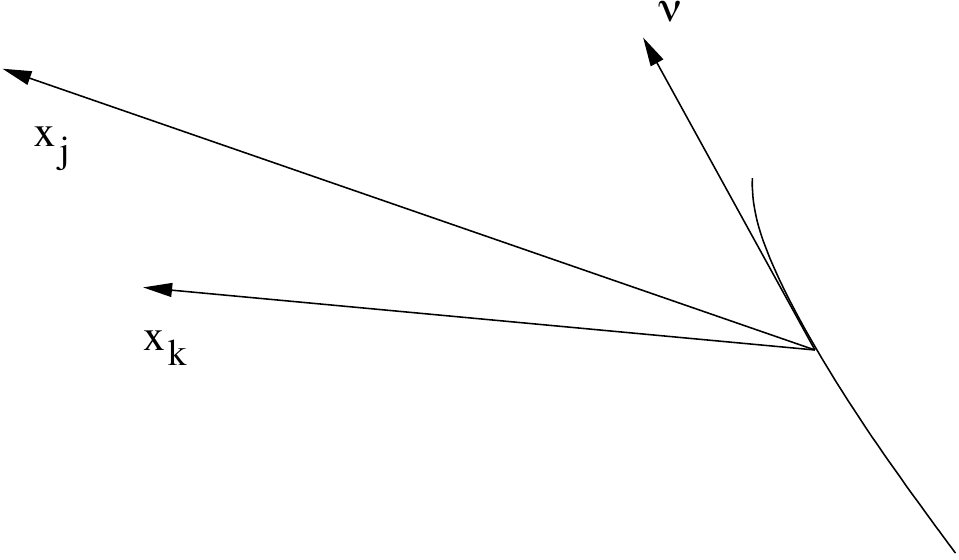}
\caption{Image of $\{(0,t)\}$ in geodesic normal coordinates about $(0,0)$}
\end{figure*}

\newsection{Local restrictions of eigenfunctions to non-smoothly closed geodesics}

We have shown above that if $\{e_{\lambda_{j_k}}\}_{k=1}^\infty$ is a sequence
of $L^2$-normalized eigenfunctions satisfying
\begin{equation}\label{3.1}\limsup_{k\to\infty}\sup_{\gamma\in\Pi}\lambda_{j_k}^{-1/2}\int_{\gamma}
|e_{\lambda_{j_k}}|^2 \, ds =0,
\end{equation}
then $\lambda_{j_k}^{-\delta(p)}\|e_{\lambda_{j_k}}\|_{L^p(M)}=0$, $2<p<6$. While it seems difficult
to determine when this  holds, one can show the following.

\begin{proposition}\label{prop3}  Suppose that $\gamma\in \Pi$ is not contained in a smoothly
closed geodesic.  Then if $\{e_{\lambda_j}\}$ is the full sequence of $L^2$-normalized
eigenfunctions, we have
\begin{equation}\label{3.2}
\limsup_{j\to\infty}\lambda_j^{-1/2}\int_{\gamma}|e_{\lambda_j}|^2 \, ds =0.
\end{equation}
\end{proposition}

In proving this proposition we may assume, after possible multiplying the metric by a constant,
that the injectivity radius is  more than 10.  This will allow us to write down Fourier integral operators 
representing the solution of the wave equation up to times $|t|\le 10$.  More important, though,
is that we shall use an observation of Tataru \cite{tataru} that the map from Cauchy data
to the solution of the wave equation restricted to $\gamma\times {\Bbb R}$ is a Fourier
integral operator with a one-sided fold.  Using this fact and the 
standard
method of long-time averages (see e.g. \cite{dg}, \cite{ivrii}, \cite{soggezelditch}, \cite{stz}),
we shall be able to prove Proposition~\ref{prop3}.

To set up our proof, let us choose Fermi normal coordinates about $\gamma$ so that, in these
coordinates, $\gamma$ becomes $\{(s,0); \, 0\le s\le 1\}$.  Note that in these coordinates the
metric takes the form $g_{11}(x)dx_1^2+dx_2^2$.  As a consequence if $p(x,\xi)=\sqrt{\sum g^{jk}(x)\xi_j\xi_k}$ is the principal symbol of $P=\sqrt{-\Delta_g}$ then $p((s,0),\xi)=\sqrt{g_{11}((s,0))\xi_1^2+\xi_2^2}$ is an even function of $\xi_2$.

To proceed, let us fix a real-valued function $\chi\in {\mathcal S}({\mathcal R})$ with $\chi(0)=1$ and
$\hat \chi(t)=0$, $|t|>1/2$.  Then if $e_\lambda$ is an eigenfunction with eigenvalue $\lambda$ it follows that $\chi(N(P-\lambda))e_\lambda=e_\lambda$.  Thus, in order to prove \eqref{3.2}, it 
would suffice to prove that given $\lambda,N\ge1$
\begin{equation}\label{g.2} \big\| \chi(N(P-\lambda))f\bigr\|_{L^2(\gamma)}
\le CN^{-1/2}\lambda^{1/4}\|f\|_{L^2(M)}+C_N\|f\|_{L^2(M)}.
\end{equation}
Note that
\begin{equation}\label{g.3}\chi(N(P-\lambda))f(x)=N^{-1}\int \hat \chi(t/N) e^{-it\lambda}\Bigl(e^{itP}f\Bigr)(x)\, dt,
\end{equation}
and because of the support properties of the $\hat \chi$ the integrand vanishes when $|t|\ge N/2$.

The operator 
$$f\to \Bigl(e^{itP}f\Bigr)(x)$$
is a Fourier operator with canonical relation
$$\{\, (x,t,\xi,\tau; y,\eta); \, \Phi_t(x,\xi)=(y,\eta), \pm \tau=p(x,\xi)\, \},$$
with $\Phi_t: \, T^*M\to T^*M$ being geodesic flow on the cotangent bundle and $p(x,\xi)$, as above,
being the principal symbol of $\sqrt{-\Delta_g}$.  Given that we want to restrict the operator in \eqref{g.3} to $\gamma=(s,0)$, $0\le s\le 1$, we really need to also focus on the the Fourier integral operator
$$f\to \Bigl(e^{itP}f\Bigr)(s,0).$$
Given the above, its canonical relation is 
$${\mathcal C}=\bigl\{\, \bigl( \, \Pi_{\gamma\times {\mathbb R}}(x,t,\xi,\tau; y,\eta\bigr)\in T^*(\gamma\times {\mathbb R})
\times T^*M; \, \Phi_t(x,\xi)=(y,\eta), \, \pm \tau=p(x_1,0,\xi)\, \bigr\},$$
with $\Pi_{\gamma\times {\mathbb R}}$ being the projection map from $T^*(M\times {\mathbb R})$
to $T^*(\gamma\times {\mathbb R})$.  Note that the projection from the latter canonical relation
to $T^*(\gamma\times{\Bbb R})$ is the map
$$(s,t,\xi)\to (s,t,\xi_1,p((s,0),\xi)),$$
which has a fold singularity when $\xi_2=0$ but has surjective differential away from this set
(given the aforementioned properties of $p$).

Because of this, given the explicit formula in Fermi coordinates, if we choose $\psi\in C^\infty_0(M)$ equal to one on $\gamma$ and $\alpha
\in C^\infty_0({\mathbb R})$ satisfying $\alpha=1$ on $[-1/2,1/2]$ but $\alpha(\tau)=0$, $|\tau|\ge1$,
then 
$$b_\varepsilon(x,\xi)=\psi(x)\alpha(\xi_2/\varepsilon|\xi|)$$
equals one on a conic neighborhood of the set that projects onto the set where the left
projection of ${\mathcal C}$ has a folding singularity.  This means that
$$B_\varepsilon(x,\xi)=\psi(x)\bigl(1-\alpha(\xi_2/\varepsilon|\xi|)\bigr)$$
has symbol vanishing in a conic neighborhood of this set and consequently the map
$$f\to \Bigl(B_\varepsilon \circ e^{itP}f\Bigr)((s,0)), \, \, \, 0\le s\le 1$$
is a nondegenerate Fourier integral operator of order zero.  Therefore, H\"ormander's theorem
\cite{hormander2} about the $L^2$ boundedness of Fourier integral operators yields
$$\int_{-N}^N \int_0^1 \, \left| \, \Bigl(B_\varepsilon \circ e^{itP}f\Bigr)(s,0)\, \right|^2 ds dt
\le C_{N,B_\varepsilon}\|f\|^2_{L^2(M)}.
$$
Therefore, an application of Schwarz's inequality yields
$$\|\chi^{N,B_\varepsilon}_\lambda f\|_{L^2(\gamma)}\le C'_{N,B_\varepsilon}\|f\|_{L^2(M)},
$$
if 
$$\chi^{N,B_\varepsilon}_\lambda f = B_\varepsilon\circ \chi(N(P-\lambda))f
=N^{-1}\int \hat \chi(t/N)e^{-it\lambda}\Bigl(B_\varepsilon\circ e^{itP}\Bigr)f dt.$$

Therefore if we similarly define $\chi^{N,b_\varepsilon}_\lambda f =b_\varepsilon\circ
\chi(N(P-\lambda))f$, then $\chi^{N,B_\varepsilon}_\lambda f+\chi^{N,b_\varepsilon}_\lambda f
=\psi \chi(N(P-\lambda))f$ and since $\psi=1$ on $\gamma$, the proof of \eqref{g.2}
would be complete if we could show that if $\varepsilon>0$ is small enough (depending
on $N$) then for $\lambda\ge1$ we have for a constant $C$ independent of $\varepsilon,
N$ and $\lambda\ge1$
\begin{equation}\label{g.4}
\|\chi^{N,b_\varepsilon}_\lambda f\|_{L^2(\gamma)}\le CN^{-1/2}\lambda^{1/4}\|f\|_{L^2(M)}
+C_{N,b_\varepsilon}\|f\|_{L^2(M)}.
\end{equation}
In addition to taking $\varepsilon>0$ to be small, we shall also take the support of $\psi$ about
$\gamma$ to be small.

It is in proving \eqref{g.4} of course where we shall use our assumption that $\gamma$ is 
not part of a smoothly closed geodesic.  A consequence of this is that, given
fixed $N$,  if $\varepsilon$ and the support of $\psi$ are small enough then
\begin{equation}\label{g.5}
b_\varepsilon(y,\eta)=0 \, \, \text{whenever } \, (y,\eta)=\Phi_t(x,\xi), 
\quad
(x,\xi)\in \text{supp }b_\varepsilon, \, \, 2\le |t|\le N.
\end{equation}
In what follows, we shall assume that $\varepsilon$ and $\psi$ have been chosen so that
this is the case.  The point here is that if $\gamma(s)$, $s\in {\mathbb R}$, is the geodesic starting at $(0,0)$
and containing $\{\gamma(s)=(s,0); \, 0\le s\le 1\}$, points on the curve $\gamma(s)$, $|s|\le N+1$
might intersect $\gamma$, but the intersection must be transverse as $s\to \gamma(s)$ is 
not a smoothly closed geodesic.  Then if $\varepsilon$ is chosen to be a small multiple of the
smallest angle of intersection and if $\psi$ has small enough support about $\gamma$, then
we get \eqref{g.5}.  Using the canonical relation for $e^{itP}$, we can deduce from this that
\begin{equation}\label{g.66}
b_\varepsilon \circ e^{itP}\circ b^*_\varepsilon \quad \text{is a smoothing operator when} \, 
2\le |t|\le N+1,
\end{equation}
i.e., for such times this operator's kernel is smooth.

Let $T$ be the operator $\chi^{N,b_\varepsilon}_\lambda f|_{\gamma},$ i.e., the
truncated approximate spectral projection operator restricted to $\gamma$.  Our goal
is to show \eqref{g.4} which says that
$$\|T\|_{L^2(M)\to L^2(\gamma)}\le CN^{-1/2}\lambda^{1/4}+C_{N,b_\varepsilon}.$$
This is equivalent to saying that the dual operator 
$T^*: L^2(\gamma)\to L^2(M)$ with the same norm, and since
$$\|T^*g\|^2_{L^2(M)}=\int_M T^*g\,  \overline{T^*g}dx=\int_\gamma TT^*g\, \, \overline{g}ds
\le \|TT^*g\|_{L^2(\gamma)}\|g\|_{L^2(\gamma)},$$
we would be done if we could show that
\begin{equation}\label{g.6}
\|TT^*g\|_{L^2(\gamma)}\le \Bigl(\, CN^{-1}\lambda^{1/2}+C_{N,b_\varepsilon} \, \Bigr)\|g\|_{L^2(\gamma)}.
\end{equation}
But the kernel of $TT^*$ is $K(\gamma(s),\gamma(s'))$, where $K(x,y)$, $x,y\in M$ is the
kernel of the operator $b_\varepsilon \circ \rho(N(P-\lambda))\circ b_\varepsilon^*$ with
$\rho(\tau)=(\chi(\tau))^2$ being the square of $\chi$.  Its Fourier transform, $\hat \rho$,
is the convolution of $\hat \chi$ with itself, and thus $\hat \rho(t)=0$, $|t|\ge 1$.  Consequently,
we can write
\begin{equation}\label{g.7}
b_\varepsilon \circ \rho(N(P-\lambda))\circ b^*_\varepsilon = 
N^{-1}\int \hat\rho(t/N)e^{-it\lambda}\Bigl(b_\varepsilon \circ e^{itP}\circ b^*_\varepsilon 
\Bigr) \, dt.
\end{equation}
Thus, if $\alpha\in C^\infty_0({\mathbb R})$ is as above, then by \eqref{g.5} and \eqref{g.66}, the difference of the kernel
of the operator in \eqref{g.7} and the kernel of the operator given by
\begin{equation}\label{g.8}
N^{-1}\int \alpha(t/10) \hat \rho(t/N) e^{-it\lambda}\Bigl( \, b_\varepsilon\circ
e^{itP}\circ b^*_\varepsilon \, \Bigr) \, dt
\end{equation}
is $O(\lambda^{-J})$ for any $J$.  Thus, if we restrict the kernel of the difference to $\gamma\times
\gamma$, it contributes a portion of $TT^*$ that maps $L^2(\gamma)\to L^2(\gamma)$ 
with norm $\le C_{N,b_\varepsilon}$.

To finish, we need to estimate the remaining piece, which has the kernel of the operator
in \eqref{g.8} restricted to $\gamma\times \gamma$.  Since we are assuming that the injectivity
radius of $M$ is 10 or more one can use the Hadamard parametrix for the wave
equation and standard stationary phase arguments (similar to ones in \cite{soggebook}, Chapter 5, or the proof of Lemma~4.1 in \cite{burq}) to see that the kernel $K(x,y)$ of the operator in \eqref{g.8} satisfies 
$$|K(x,y)|\le CN^{-1}\lambda^{1/2}\bigl(d_g(x,y)\bigr)^{-1/2}+C_{b_\varepsilon}.$$
The first term comes from the main term in the stationary phase expansion for the kernel and
the other one is the resulting remainder term in the one-term expansion.  Since this kernel
restricted to $\gamma\times \gamma$ gives rise to an integral operator satisfying
the estimates in \eqref{g.6}, the proof is complete. \qed

\newsection{Further questions}

While as we explained before the condition that for the $L^2$-normalized eigenfunctions
$$\limsup_{j\to \infty}\sup_{\gamma\in \Pi}\lambda^{-1/2}_j\int_\gamma |e_{\lambda_j}|^2\, ds=0
$$
is a natural one to quantify non-concentration, it would be interesting to formulate a geometric condition involving the
long-time dynamics of the geodesic flow that would imply it and its equivalent version
that $\lambda_j^{-\delta(p)}\|e_{\lambda_j}\|_p\to 0$, $2<p<6$.  Presumably if $\gamma\in \Pi$
and 
\begin{equation}\label{4.1}\limsup_{j\to \infty}\lambda^{-1/2}_j\int_\gamma |e_{\lambda_j}|^2 ds>0,\end{equation}
 then $\gamma$
would have to be part of a stable smoothly closed geodesic, and not just a closed geodesic as we showed above.  
Toth and Zeldtich  made a similar conjecture to this in \cite{tothzelditch}, saying that, in $n$-dimensions, if $\gamma$ is a closed stable geodesic then one should be able to find a sequence of eigenfunctions on which sup-norms are blowing up like $\lambda^{(n-1)/2}$.  In \cite{ralstonb}, \cite{ralston}, it was shown that there is a sequence of quasimodes blowing up at this rate.

It would also be interesting to
formulate a condition that would ensure that $\|e_\lambda\|_{L^6(M)}=o(\lambda^{\delta(6)})=o(\lambda^{1/6})$, for $L^2$-normalized eigenfunctions.  Presumably, such a condition would have to involve both ones like those in the present paper and conditions of the type in  \cite{stz}, \cite{soggezelditch}.  Since $L^6$
is an endpoint for \eqref{1} one expects that one would need a condition that both guarantees
that $L^p$ bounds for $2<p<6$ and $p>6$ be small.  Formally, the proof of Theorem~\ref{theorem1} suggests that $L^4$-norms over geodesics might be relevant for the problem
of determining when the $L^6(M)$ norms of eigenfunctions are small.  This is interesting
because the $L^4$-norm is the unique $L^p$-norm taken over geodesics that captures both
the concentration of the highest weight spherical harmonics on geodesics and the concentration
of zonal functions at points.  Indeed, the highest weight spherical harmonics saturate these
norms for $2\le p\le 4$, while the zonal functions saturate them for $p\ge4$ (see \cite{burq}).

Also, it would be interesting to see whether the results here generalize to the case of two-dimensional
compact manifolds with boundary.  Recently, Smith and the author \cite{ss} were able to obtain sharp eigenfunction estimates
in this case.  In this case, the critical estimate was an $L^8$ one.  So the results here suggest
that size estimates for the Kakeya-Nikodym maximal operator associated with broken unit geodesics
and applied to squares of eigenfunctions could be relevant for improving the bounds in \cite{ss}, which are known to be sharp in the case of
the disk (see \cite{grieser}).  An observation
of Grieser \cite{grieser} involving the Rayleigh whispering gallery modes suggests that
in order to obtain a variant of Corollary~\ref{corollary1} for compact
domains one would have to consider $L^2$-norms
over $\lambda^{-2/3}_j$-neighborhoods of broken geodesics.  Smith and the author~\cite{ss0} also showed that for compact manifolds with geodesically concave boundary one has better
estimates than one does for compact domains in ${\mathbb R}^n$.  For example, when $n=2$
\eqref{1} holds.  Based on this and the better behavior of the geodesic flow, it seems
reasonable that the analog of Corollary~\ref{corollary1} might hold (with the same scales)
in this setting.

Finally, as mentioned before it would be interesting to see to what extent the results for the boundaryless case
extend to higher dimensions.  The arguments given here and in \cite{bourgainef}, though,
rely very heavily on special features of the two-dimensional case.  

\bigskip
\noindent{\bf Acknowlegements:}
It is a real pleasure to thank J. Bourgain for sharing an early version of his paper \cite{bourgainef} and to also thank W. Minicozzi  for helpful conversations and for going over a key step in the proof.  The author would also like to express his gratitude to J.~Toth and S. Zeldtich for helpful discussions and suggestions.


\begin{thebibliography}{MA}
\bibitem{ralstonb} V. M. Babi\v c and V. F. Lazutkin, {\em The eigenfunctions which are
concentrated near a closed geodesic}, (Russian)  Zap. Nau\v cn. Sem. Leningrad. Otdel. Mat. Inst. Steklov, (LOMI)  {\bf 9}  1968 15--63.
\bibitem{bourgain1} J. Bourgain, {\em Besicovitch type maximal operators and applications to Fourier analysis},  Geom. Funct. Anal.  {\bf 1}  (1991),   147--187.
\bibitem{bourgain2} J. Bourgain, {\em Some new estimates on oscillatory integrals}, Annals of Math.
Studies, {\bf 42} (1995), 83--112.
\bibitem{bourgain3} J. Bourgain, {\em $L\sp p$-estimates for oscillatory integrals in several variables},
Geom. Funct. Anal. {\bf 1} (1991), 321--374. 
\bibitem{bourgainef} J. Bourgain, {\em Geodesic restrictions and $L^p$-estimates
for eigenfunctions of Riemannian surfaces}, Linear and Complex Analysis: Dedicated to V. P. Havin on the Occasion of His 75th Birthday, American Math. Soc. Transl., Advances in the Mathematical Sciences (2009), 27-35.
\bibitem{burq} N. Burq, P. G\'erard and N. Tzvetkov, {\em Restriction of the Laplace-Beltrami eigenfunctions to submanifolds}, Duke Math. J.  {\bf 138}
(2007), 445--486.
\bibitem{carsj} L. Carleson and P. Sj\"olin, {\em Oscillatory integrals and a multiplier 
problem for the disc}, Sudia Math. {\bf 44} (1972), 287--299.
\bibitem{cdv} Y. Colin de Verdi\`ere, {\em Semi-classical measures and entropy [after Nalini Anantharaman and St\'ephane Nonnenmacher]},  (English summary)
S\'eminaire Bourbaki. Vol. 2006/2007.
AstŽrisque No. 317 (2008), Exp. No. 978, ix, 393--414.
\bibitem{cordoba} A. C\'ordoba, {\em A note on Bochner-Riesz operators}, Duke Math. J. {\bf 46} (1979), 505--511.
\bibitem{dg} J. J. Duistermaat and V. W. Guillemin, {\em The spectrum of positive elliptic operators and periodic bicharacteristics},  Invent. Math.  {\bf 29}  (1975), 39--79.
\bibitem{fefferman} C. Fefferman, {\em A note on spherical summation operators}, Israel
J. Math. {\bf 15} (1973), 44-52.
\bibitem{greenleafseeger} A. Greenleaf and A. Seeger, {\em Fourier integrals with fold
singularities}, J. Reine Angew. Math. {\bf 455} (1994), 35--56.
\bibitem{grieser} D. Grieser, {\em
$L^p$ Bounds for Eigenfunctions and Spectral Projections of the Laplacian Near Concave Boundaries},  Ph. D. thesis, University of California: Los Angeles, 1992.
\bibitem{hormander2} L. H\"ormander, {\em Fourier integral operators. I},  Acta Math.  {\bf 127}  (1971), 79--183.
\bibitem{hormander} L. H\"ormander, {\em Oscillatory integrals and multipliers on $FL^p$},
Ark. Math. II (1973), 1--11.
\bibitem{ivrii} V. Ivrii, {\em The second term of the spectral asymptotics for a Laplace-Beltrami operator on manifolds with boundary} (Russian)  Funktsional. Anal. i Prilozhen.  {\bf 14}  (1980), 25--34.
\bibitem{mins} W. P. Minicozzi and C. D. Sogge, {\em Negative results for Nikodym maximal functions and related oscillatory integrals in curved space},  Math. Res. Lett.  {\em 4}  (1997),  221--237.
\bibitem{mss} G. Mockenhaupt, A. Seeger, C. D. Sogge, {\em Local smoothing of Fourier
integral operators and Carleson-Sj\"olin estimates}, J. Amer. Math. Soc. {\bf 6} (1993), 65--130.
\bibitem{ralston} J. V. Ralston, {\em On the construction of quasimodes associated with stable periodic orbits},  Comm. Math. Phys.  {\bf 51}  (1976), 219--242. 
\bibitem{rez} A. Reznikov, {\em Norms of geodesic restrictions for eigenfunctions on hyperbolic
surfaces and representation theory}, arXiv:math.AP/0403437.
\bibitem{sarnak} Z. Rudnick and P. Sarnak, {\em The behaviour of eigenstates of arithmetic
hyperbolic manifolds}, Comm. Math. Phys {\bf 161} (1994), 195--213.
\bibitem{schnirelman} A. Schnirelman: {\em Ergodic properties of eigenfunctions}, Usp. Math. Nauk. {\bf 29}, (1974), 181--182.
\bibitem{ss0} H. Smith and C. D. Sogge, {\em On the critical semilinear wave equation outside convex obstacles},
J. Amer. Math. Soc. {\bf 8} (1995), 879--916. 
\bibitem{ss} H. Smith and C. D. Sogge, {\em On the $L\sp p$ norm of spectral clusters for compact manifolds with boundary}, Acta Math. {\bf 198} (2007), 107--153.
\bibitem{sph} C. D. Sogge, {\em Oscillatory integrals and spherical harmonics}  Duke Math. J.  
{\bf 53}  (1986),  43--65.
\bibitem{soggeest} C. D. Sogge, {\em Concerning the $L^p$ norm of spectral clusters for second-order elliptic operators on compact manifolds}, J. Funct. Anal. {\bf 77} (1988), 123--138.
\bibitem{soggebook} C. D. Sogge, {\em Fourier integrals in classical analysis}, Cambridge Tracts in Math., Cambridge Univ. Press, Cambridge, 1993.
\bibitem{sonick} C. D. Sogge, {\em Concerning Nikodym-type sets in $3$-dimensional curved spaces},  J. Amer. Math. Soc.  {\bf 12}  (1999), 1--31.
\bibitem{stz} C. D. Sogge, J. Toth and S. Zelditch, {\em About the blowup of quasimodes
on Riemannian manifolds}, to appear, J. Geom. Anal.
\bibitem{soggezelditch} C. D. Sogge and S. Zelditch, {\em Riemannian manifolds with maximal
eigenfunction growth}, Duke Math. J. {\bf 114} (2002), 387--437.
\bibitem{stein} E. M. Stein, {\em Oscillatory integrals in Fourier analysis}, Beijing Lectures
in Harmonic Analysis, Princeton Univ. Press, Princeton, NJ, 1986, pp. 307--356.
\bibitem{tataru} D. Tataru, {\em On the regularity of boundary traces for the wave equation}, 
Ann. Scuola Norm. Sup. Pisa Cl. Sci.   {\bf 26}  (1998),   185--206. 
\bibitem{toth} J. Toth, {\em $L^2$-restriction bounds for eigenfunctions along curves in
the quantum completely integrable case}, Comm. Math. Phys. {\bf 288} (2009), 379--401.
\bibitem{tothzelditch} J. Toth and S. Zelditch, {\em $L^p$ norms of eigenfunctions in the completely integrable case}, Ann. Henri Poincar\'e {\bf 4} (2003), 343--368.
\bibitem{zelditch} S. Zelditch, {\em Uniform distribution of eigenfunctions on compact hyperbolic surfaces},  Duke Math. J.  {\bf 55}  (1987),  919--941. 
\bibitem{zelind} S. Zelditch, {\em Mean Lindel\"of hypothesis and equidistribution of cusp forms and Eisenstein series},
J. Funct. Anal. {\bf 97} (1991), 1--49. 
\bibitem{zygmund} A. Zygmund, {\em On Fourier coefficients and transforms of 
two variables}, Studia Math. {\bf 50} (1974), 189--201.


\end{thebibliography}
\end{document}